\documentclass[letterpaper,10pt,conference]{ieeeconf}

\IEEEoverridecommandlockouts
\usepackage{mathtools,amssymb}
\usepackage{graphicx}
\graphicspath{{figures/}}
\usepackage{cite}
\newtheorem{proposition}{Proposition}
\newtheorem{lemma}{Lemma}

\title{\LARGE\bfseries
 On the minimum
energy required to trigger a spike
}

\author{Guanchun Tong$^{1}$ and Rodolphe Sepulchre$^{2}$%
\thanks{This work was supported by the European Research Council under the
Advanced ERC Grant Agreement SpikyControl (project number: 101054323).
The work of Guanchun Tong was also supported by Flemish Fund for Scientific
Research (FWO) under a PhD Fellowship fundamental research (project number:
1191125N).}%
\thanks{$^{1}$Guanchun Tong is with the Department of Electrical Engineering
(ESAT), KU Leuven, Belgium.
\hfill\break\texttt{guanchun.tong@esat.kuleuven.be}}%
\thanks{$^{2}$Rodolphe Sepulchre is with the Department of Electrical Engineering
(ESAT), KU Leuven, Belgium, and also with the Department of
Engineering, Control Group, University of Cambridge, UK.
\hfill\break\texttt{rodolphe.sepulchre@kuleuven.be}}%
}

\begin{document}

\maketitle
\thispagestyle{empty}
\pagestyle{empty}

\begin{abstract}

We formulate and numerically study the problem of finding the minimum energy required to trigger a 
spike in a three-dimensional conductance-based excitable model. 
The variational formulation seeks a state trajectory that reaches the boundary of 
the spiking region with the least \textit{peak} cumulative work supplied through the membrane port.
We approximate voltage histories by sums of real exponentials, reducing the problem
to a finite-dimensional optimization. 
The best-found voltage waveform consists of a hyperpolarization 
phase followed by a depolarization phase, which seems to match the 
experimentally determined optimal stimulus in squid giant axons.

\end{abstract}

\section{Introduction}
\label{sec:introduction}

In 1744, Leonhard Euler wrote that ``Nothing at all happens in the world in which 
some principle—whether of the maximum or the minimum—does not manifest itself''~\cite[p.~245]{Euler1744}. 
This conviction captures
the unifying power of variational principles: physical behavior can be
characterized through a single scalar functional of its trajectories. 
From classical mechanics to quantum theory and general relativity, variational
formulations underpin some of our most fundamental descriptions of physical laws.

Our objective is to apply this principle to neuronal excitability, a nonlinear phenomenon
characterized by its all-or-none response to a stimulus--the threshold property. 
The mathematical study of excitability
started in the seminal work of Hodgkin and Huxley on the quantitative modeling
 of the electrical activity in the squid giant axon in 1952~\cite{HodgkinHuxley1952}. 
 More than seventy
years later, their theory remains arguably still the \textit{only} triumph
in mathematical neuroscience. Excitability is also observed in cardiac
cells, muscle cells, and other cells in which electrical
signaling is important. The ubiquity of excitability
calls for a deeper understanding of its underlying principles, ideally 
through a variational principle.
The Hodgkin--Huxley theory, however, does not itself provide a variational
formulation.

Two facts complicate such a formulation: excitable systems
exchange energy with their environment while their active elements
can also supply energy internally. The circuit representation of
Hodgkin--Huxley nevertheless makes it amenable to an energy-based analysis using
open-system theory~\cite{SepulchreDrionFranci2018}.
In our earlier work~\cite{SepulchreTong2025}, we proposed
characterizing the threshold by a local maximum of a voltage-dependent value 
function, which is closely related to the notion of required supply, a central 
concept of dissipativity theory~\cite{Willems1972,VanDerSchaft2021}. 

However, our previous characterization focused on the terminal voltage, whereas 
required supply is a function of the entire state.
This article continues that required supply viewpoint by studying the energy required to 
reach the spiking region in the full state space, linking the energy criterion
 to the global geometry of spike initiation.

Understanding which features of a stimulus lead to a spike is a
longstanding question in neuronal computation~\cite{MainenSejnowski1995,Koch1999}.
The metabolic cost of neural processing motivates an energy-based
view of this question~\cite{NivenLaughlin2008}, while optimal-control
methods have been used to design energy-efficient excitation
currents~\cite{JezernikMorari2005}. Variational optimization of root
mean square (RMS) current has also been combined with experiments
on squid giant axons~\cite{ForgerPaydarfarClay2011}. 

Moreover, an energy-based analysis of excitable systems is relevant
to neuromorphic engineering, where energy efficiency is a central
motivation. The tradeoff between energy consumption and computational
performance is a standard part of quantitative computer
architecture~\cite{HennessyPatterson2017}. In analog spiking neurons,
energy dissipation has also been related to spike-timing reliability~\cite{VanBrandtEtAl2026}, and therefore
related to the post-inhibitory rebound mechanism~\cite{Sejnowski2026}. Our work provides
a variational framework for studying the energy
required to trigger a spike, where the importance of inhibition is revealed 
in the computational results.

We begin in Section~\ref{sec:structure} by introducing a three-dimensional 
model of an excitable system, the recently proposed excitatory--inhibitory 
(EI) model~\cite{BurgerSepulchre2026}. We choose this model for three reasons. 
First, it retains the conductance-based fast--slow mixed feedback structure of
Hodgkin--Huxley, providing a physical circuit model where the open-system nature 
is represented by the port and active components are simply the batteries.
Second, it localizes the nonlinearities of the gating dynamics to a single 
static nonlinearity so that the relationship between the membrane voltage and 
the internal dynamic variables is linear. Third, its three-dimensional state 
space allows us to apply nullcline analysis to the limit of direct 
visualization, providing geometric insight into the numerically computed 
zero-input spiking region that serves as an oracle for the variational problem.

In Section~\ref{sec:variational-problem}, we formulate the variational problem
that addresses the two challenges of open systems and active components. To 
mitigate the effect of internal batteries, we introduce a novel objective 
function that measures the maximum cumulative work supplied over the course of
preparation. Subsequent extraction of energy from the active components
cannot reduce a previously recorded peak. To address the system's open nature, 
 we observe that the input current influences the internal gating variables
 only through the membrane voltage. We can therefore treat the
 membrane voltage history as the free variable to be optimized over, 
 effectively inverting the system. The external current can then be 
 reconstructed from the optimized voltage trajectory.

As this is a highly nonconvex and infinite-dimensional problem, 
 we have to turn to approximation and numerical optimization. Section~\ref{sec:real-exponential-approximation} introduces an
approximation of voltage signals by a sum of real exponentials, reducing the
infinite-dimensional optimization problem to a finite-dimensional one.
 However, the approximations should not compromise the original question we want to 
 answer, that is, \textit{What is the minimum energy required to trigger a spike?}
Crucially, for each candidate voltage trajectory, there are two consequences we are
interested in: the continuous one of energy considerations and the discrete one 
of whether the trajectory leads to a spike. We prove that this approximation will not distort 
the discrete aspect of the question: three exponentials can reach
any terminal point in the three-dimensional state space, therefore covering 
all states in the spiking region.

In Section~\ref{sec:numerical-results}, we present the numerical results of the
optimization problem. The computed voltage waveform consists of a hyperpolarization 
phase followed by a depolarization phase, which echoes the mixed feedback 
architecture of the circuit and appears to be biologically plausible. The 
continuous aspect of the above question is also tested numerically: we find that 
increasing the number of exponentials does not further reduce the optimal value 
of the energy objectives.
\section{Structure and Geometry of the EI Model}
\label{sec:structure}

\begin{figure*}[t]
    \centering
    \includegraphics[width=0.9\textwidth]{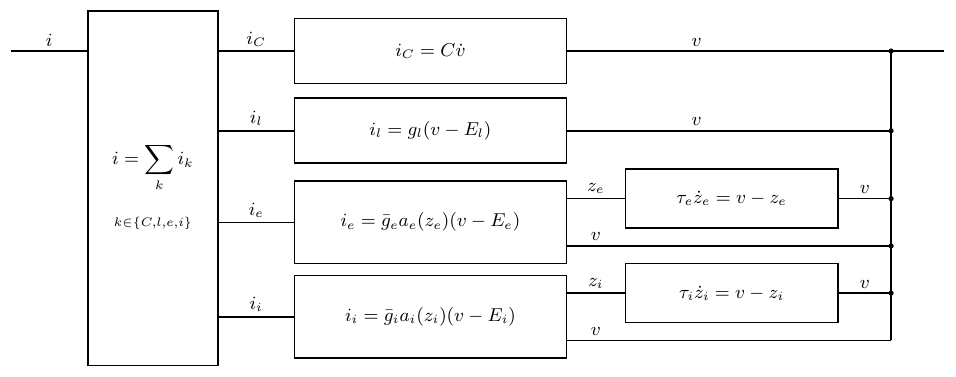}
    \caption{Behavioral constraint graph of the EI model.
    The left block enforces $i=i_C+i_l+i_e+i_i$.
    Lines connect the blocks, and dots distribute the common
    voltage signal. Each nonlinear current map receives both $v$ directly
    and its filtered state $z_k$. Branch currents use the outward-current
    sign convention. The activation functions are
    $a_k(z_k)=[1+\exp(-\kappa(z_k-\theta_k))]^{-1}$ for $k\in\{e,i\}$,
    with $\tau_e<\tau_i$.}
    \label{fig:ei-signal-flow}
\end{figure*}

We consider the temporal excitatory--inhibitory (EI) model described
in~\cite{BurgerSepulchre2026}. Its behavioral constraint graph is shown in
Fig.~\ref{fig:ei-signal-flow}. The applied current $i$ and the returning
branch currents determine the capacitor current
$i_C=i-i_e-i_i-i_l$, which drives the membrane voltage $v$
through $C\dot v=i_C$. Each nonlinear
branch current depends on both its filtered state and the instantaneous
voltage; the leak current depends directly on voltage.

\subsection{Nullcline Analysis}

The nullcline surfaces of the zero-input EI model are defined by the equations
$\dot v=0$, $\dot z_e=0$, and $\dot z_i=0$.
Write $G_k(z)=\bar g_k\sigma_k(z)$, where $\sigma_k=a_k$ is the sigmoid
in Fig.~\ref{fig:ei-signal-flow}, for $k\in\{e,i\}$.
The voltage nullcline surface is given by the equation
\begin{equation}
  -\bar{g}_e \sigma_e(z_e) (v - E_e) - \bar{g}_i \sigma_i(z_i) (v - E_i)
  - g_l (v - E_l) = 0,
\end{equation}
which can be rearranged to express $v$ explicitly as a function of $z_e$
and $z_i$, thanks to the conductance-based formulation of the model. The
resulting formula is
\begin{equation}
  v= N(z_e, z_i)=  \frac{\bar{g}_e \sigma_e(z_e) E_e + \bar{g}_i \sigma_i(z_i)
  E_i + g_l E_l}{\bar{g}_e \sigma_e(z_e) + \bar{g}_i \sigma_i(z_i) + g_l}.
\end{equation}

For fixed values of $z_e$ and $z_i$, the nullcline voltage $N(z_e, z_i)$ is a convex
combination of the reversal potentials $E_e$, $E_i$, and $E_l$; therefore it is
bounded between the minimum and maximum of these potentials, as illustrated in
Fig.~\ref{fig:nullcline-hysteresis}. The surface is smooth.

\subsection{Frozen-inhibition hysteresis}
\begin{figure}[t]
    \centering
    \includegraphics[width=0.9\columnwidth]{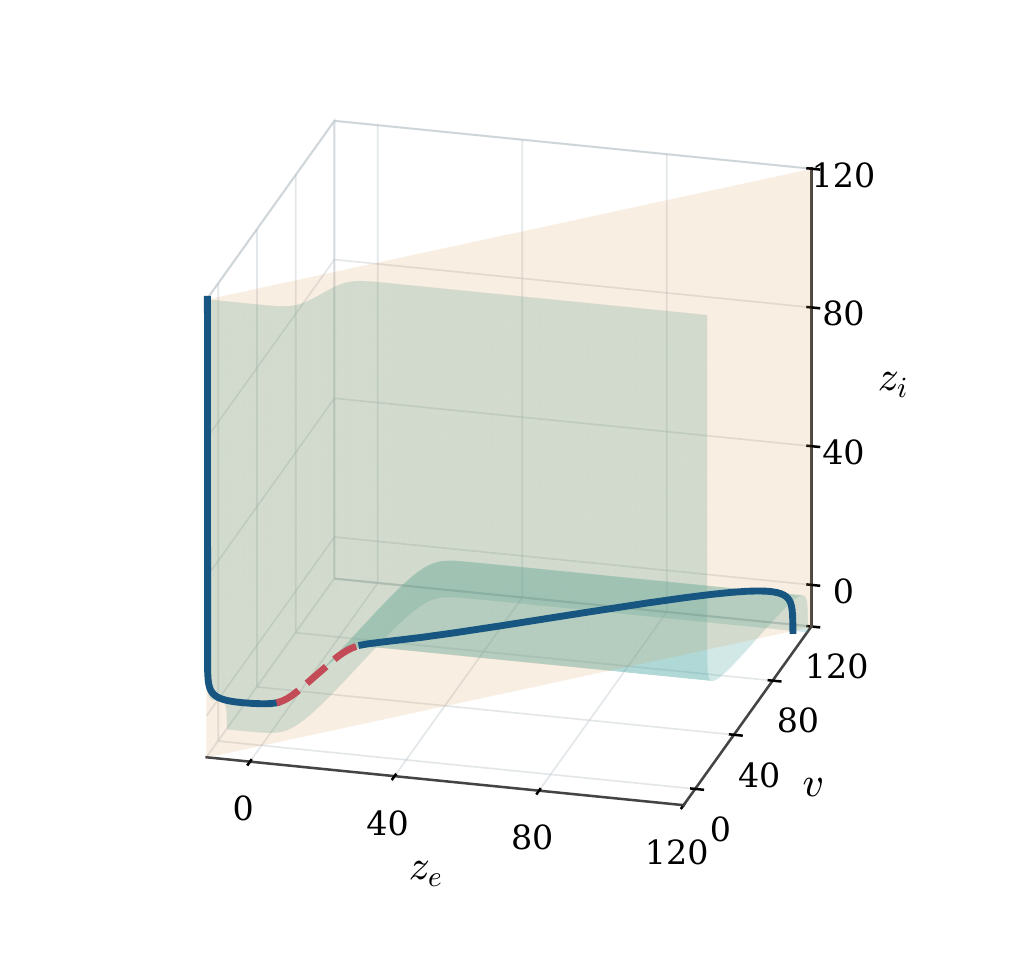}
    \caption{The voltage nullcline (teal), the plane $v=z_e$ (orange),
    and their intersection. Solid blue branches are attracting and the
    dashed red branch consists of saddle equilibria for the two-state system with $z_i$
    held fixed. The two folds produce the frozen-inhibition hysteresis
    described in the text. Parameters are those of
    Section~\ref{sec:numerical-results}.}
    \label{fig:nullcline-hysteresis}
\end{figure}

However, intersecting the voltage nullcline with the plane $v=z_e$ yields the
equilibrium curve in Fig.~\ref{fig:nullcline-hysteresis}.

\subsection{The spiking manifold}
Write $x=(v,z_e,z_i)^\top$ and let $\phi(t;x)$ denote the zero-input flow released from $x$.
 A state belongs to $\mathcal S$ if its
zero-input release reaches $v_{\mathrm{th}}$ within time $T$. That is,
\begin{equation}
\label{eq:release-spike-set}
    \mathcal S=\{x\in [E_{\mathrm{i}},E_{\mathrm{e}}]^3:\exists t\in[0,T],[\phi(t;x)]\vert_{v}\geq v_{\mathrm{th}}\}.
\end{equation}

\begin{figure}[!htbp]
    \centering
    \includegraphics[width=0.9\columnwidth]{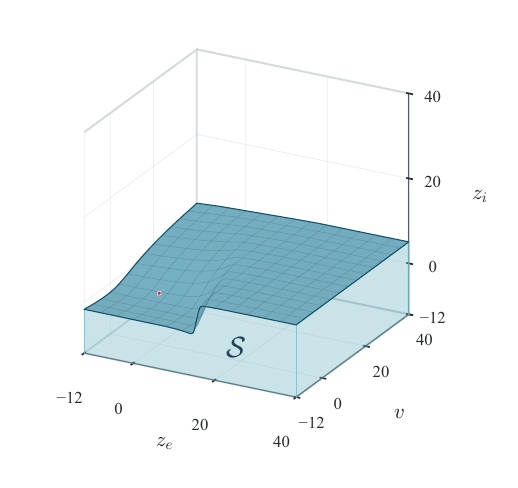}
    \caption{Numerical approximation of the zero-input spiking region
    $\mathcal S$ for $T=50\,\mathrm{ms}$ and
    $v_{\mathrm{th}}=75\,\mathrm{mV}$, restricted to initial states in
    $[-12,40]^3$. The blue sheet is the interpolated release threshold,
    and the shaded volume represents its spiking side. Pale side and
    bottom faces mark the viewing limits; the red point marks rest.
    Coordinates are in mV.}
    \label{fig:spiking-manifold}
\end{figure}

Informally, we will refer to $\mathcal S$ as the \emph{spiking manifold}, which is a three-dimensional connected region
within the cube $[E_{\mathrm{i}},E_{\mathrm{e}}]^3$.
We can visualize the spiking manifold simply by solving the zero-input dynamics for a grid of initial
conditions in the cube and checking whether the voltage reaches $v_{\mathrm{th}}$ within time $T$.
For the parameters used throughout the paper, with $T=50\,\mathrm{ms}$ and
$v_{\mathrm{th}}=75\,\mathrm{mV}$, a numerical approximation of the spiking
region is shown in Fig.~\ref{fig:spiking-manifold}.

The boundary $\partial\mathcal S$ is a two-dimensional surface that separates the spiking and non-spiking
regions.

\section{Formulation of the Variational Problem}
\label{sec:variational-problem}

Our goal is to find the minimum energy input that can trigger a spike in the EI model. 
We formulate this as an optimization problem where we minimize the energy supplied via the input port
 subject to the dynamics of the EI model and the resting and spiking conditions.

The optimization problem can be stated as follows:
\begin{equation}
\begin{aligned}
\inf \quad & \sup_{t\leq0}\int_{-\infty}^{t} i(\tau)v(\tau)\,d\tau \\
\text{subject to:} \quad & (i, v, z_e, z_i) \in \mathfrak{B}\vert_{i\in \mathcal{I}}, \\
& (v, z_e, z_i)(-\infty) = (0, 0, 0), \\
& (v, z_e, z_i)(0) \in \mathcal{S},
\end{aligned}
\end{equation}
where $\mathcal{I}$ is the set of admissible input currents, and $\mathfrak{B}\vert_{i\in \mathcal{I}}$ 
denotes the set of all possible trajectories satisfying the dynamics and subject
to the input constraints $i \in \mathcal{I}$. $\mathcal{S}$ is the set of 
terminal states whose zero-input trajectories reach $v_{\mathrm{th}}$ within time $T$, 
as defined in Section~\ref{sec:structure}. 

The key insight of the structure of the optimization problem is that the input current $i$ only 
influences the dynamics through the voltage $v$, which in turn affects the filtered gating variables $z_e$ and $z_i$.
This can be stated as the following lemma:
\begin{lemma}
Let $\mathcal{V} = \{v: (i, v, z_e, z_i) \in \mathfrak{B}\vert_{i\in \mathcal{I}}\}$. Then 
$\mathfrak{B}\vert_{i\in \mathcal{I}} = \mathfrak{B}\vert_{v \in \mathcal{V}}$.
\end{lemma}

This is left invertibility in input-output systems theory. The 
dependence of the gating dynamics on voltage alone
also suggests a deterministic analogy with conditional
independence in causality~\cite{Pearl2009}: once the voltage history is specified, 
the input current provides no additional information about 
 the dynamics of the internal gating variables. 

Therefore, we can reformulate the optimization problem in terms of the voltage signal $v$ as follows:
\begin{equation}
\begin{aligned}
\inf \quad & \sup_{t\leq0}\int_{-\infty}^{t} i(\tau)v(\tau)\,d\tau \\
\text{subject to:} \quad & (i, v, z_e, z_i) \in \mathfrak{B}\vert_{v\in \mathcal{V}}, \\
& (v, z_e, z_i)(-\infty) = (0, 0, 0), \\
& (v, z_e, z_i)(0) \in \mathcal{S}.
\end{aligned}
\end{equation}

By the KCL constraint, the input current $i$ can be expressed as a function of the 
voltage $v$ and the filtered gating variables $z_e$ and $z_i$:
\begin{equation}
i = i_C(v) + i_l(v) + i_e(v, z_e) + i_i(v, z_i),
\end{equation}
where $i_C(v) = C \dot{v}$, $i_l(v) = g_l (v - E_l)$, $i_e(v, z_e) = \bar{g}_e \sigma_e(z_e) (v - E_e)$, 
and $i_i(v, z_i) = \bar{g}_i \sigma_i(z_i) (v - E_i)$.

The optimization problem can be further simplified as follows:
\begin{equation}
\begin{aligned}
\inf_{v \in \mathcal{V}} \quad & \sup_{t\leq0}\int_{-\infty}^{t} \bigl[i_C(v) + i_l(v)\\
&\qquad {}+ i_e(v, z_e) + i_i(v, z_i)\bigr]v\,d\tau \\
\text{subject to:} \quad & \tau_e\dot{z}_e = -z_e + v, \\
& \tau_i\dot{z}_i = -z_i + v, \\
& (v, z_e, z_i)(-\infty) = (0, 0, 0), \\
& (v, z_e, z_i)(0) \in \mathcal{S},
\end{aligned}
\end{equation}
where the dynamics of $z_e$ and $z_i$ are linear and can be solved 
explicitly given the voltage signal $v$.

We could further simplify the optimization problem by noticing that our 
objective function in fact allows us to consider only those terminal states
that are on the boundary of the spiking manifold $\mathcal{S}$. 
We state this as the following lemma:

\begin{lemma}
Assume that $\mathcal S$ is closed, the resting state lies outside
$\mathcal S$, and admissible histories remain admissible under
truncation and time translation.
Then the constraint $(v, z_e, z_i)(0) \in \mathcal{S}$ can be replaced by
$(v, z_e, z_i)(0) \in \partial \mathcal{S}$ without changing the optimal 
value of the optimization problem.
\end{lemma}

\begin{proof}
    Let $v \in \mathcal{V}$ be a voltage history that takes the system from the
    resting state $(0, 0, 0)$ to a terminal state $(v, z_e, z_i)(0) \in \mathrm{int} 
    \mathcal{S}$. Because the spiking manifold $\mathcal{S}$ is separated from the 
    resting state by its boundary $\partial \mathcal{S}$, there exists a time 
    $t_0 < 0$ such that $(v, z_e, z_i)(t_0) \in \partial \mathcal{S}$. If the trajectory
    after $t_0$ cannot further increase the inner part of the objective function 
    $\sup_{t\leq 0} \int_{-\infty}^{t} i(\tau)v(\tau)\,d\tau$, then for this $v$, we can 
    simply consider the past up to the boundary time $t_0$. If the trajectory
    after $t_0$ can further increase the inner part of the objective function, then
    outer ($\inf$) part of the optimization problem will rule out this $v$ anyway, 
    since we can always find a voltage history that takes the system only to the boundary $\partial \mathcal{S}$,
    along the same trajectory and just stops.
\end{proof}

The optimization problem now reads
\begin{equation}
\begin{aligned}
\inf_{v \in \mathcal{V}} \quad & \sup_{t\leq0}\int_{-\infty}^{t} \bigl[i_C(v) + i_l(v)\\
&\qquad {}+ i_e(v, z_e) + i_i(v, z_i)\bigr]v\,d\tau \\
\text{subject to:} \quad & \tau_e\dot{z}_e = -z_e + v, \\
& \tau_i\dot{z}_i = -z_i + v, \\
& (v, z_e, z_i)(-\infty) = (0, 0, 0), \\
& (v, z_e, z_i)(0) \in \partial \mathcal{S},
\end{aligned}
\end{equation}

\section{Approximation of the Voltage Signal Space with Real Exponentials}
\label{sec:real-exponential-approximation}
We approximate the space of admissible voltage signals with real exponential 
functions, as in \cite{SepulchreTong2025}.

In fact, real exponentials are eigenfunctions for linear relaxation systems.
Let
\begin{equation}
    v(t) = \sum_{k=1}^{N} c_k e^{\alpha_k t}.
\end{equation}
The filtered gating variables become
\begin{equation}
    z_e(t) = \sum_{k=1}^{N}\frac{c_k}{1+\alpha_k \tau_e} e^{\alpha_k t}, \quad
    z_i(t) = \sum_{k=1}^{N}\frac{c_k}{1+\alpha_k \tau_i} e^{\alpha_k t}.
\end{equation}

We can now state a proposition that characterizes the number of real exponentials
needed to approximate the voltage signal such that the resulting state trajectory 
can reach the whole state space.

\begin{proposition}
\label{prop:real-exponential-approximation}
For any three distinct positive rates $\alpha_1,\alpha_2,\alpha_3$
and any prescribed terminal state $x_f$,
there exist unique real coefficients $c_1,c_2,c_3$ such that the
voltage history $v(t)=\sum_{k=1}^{3}c_k e^{\alpha_k t}$ drives the
EI model from rest as $t\to-\infty$ to $x(0)=x_f$.
\end{proposition}

\begin{proof}
    For $N=3$, these relations have the matrix form
\begingroup
\setlength{\arraycolsep}{3pt}
\begin{equation}
\label{eq:three-exponential-state-map}
\begin{bmatrix}
    v(t) \\ z_e(t) \\ z_i(t)
\end{bmatrix}
=
\begin{bmatrix}
    e^{\alpha_1 t} & e^{\alpha_2 t} & e^{\alpha_3 t} \\[6pt]
    \dfrac{e^{\alpha_1 t}}{1+\alpha_1\tau_e} &
    \dfrac{e^{\alpha_2 t}}{1+\alpha_2\tau_e} &
    \dfrac{e^{\alpha_3 t}}{1+\alpha_3\tau_e} \\[10pt]
    \dfrac{e^{\alpha_1 t}}{1+\alpha_1\tau_i} &
    \dfrac{e^{\alpha_2 t}}{1+\alpha_2\tau_i} &
    \dfrac{e^{\alpha_3 t}}{1+\alpha_3\tau_i}
\end{bmatrix}
\begin{bmatrix}
    c_1 \\ c_2 \\ c_3
\end{bmatrix}.
\end{equation}
At $t=0$, the terminal state is
\begin{equation}
\label{eq:three-exponential-terminal-map}
x_f= A_0 c,
\end{equation}
where
\begin{equation}
A_0 = \begin{bmatrix}
    1 & 1 & 1 \\[2pt]
    \dfrac{1}{1+\alpha_1\tau_e} &
    \dfrac{1}{1+\alpha_2\tau_e} &
    \dfrac{1}{1+\alpha_3\tau_e} \\[4pt]
    \dfrac{1}{1+\alpha_1\tau_i} &
    \dfrac{1}{1+\alpha_2\tau_i} &
    \dfrac{1}{1+\alpha_3\tau_i}
\end{bmatrix}, \quad
c = \begin{bmatrix}
    c_1 \\ c_2 \\ c_3
\end{bmatrix}.
\end{equation}
\endgroup
We are going to show that $A_0$ has full rank, which implies $A_0$ is invertible, and thus for any terminal state
$(v(0), z_e(0), z_i(0))$, there exists a unique set of coefficients $c_k$ such that 
$(v(0), z_e(0), z_i(0)) = A_0 c$. Assume we have a non-zero vector $r=(r_1,r_2,r_3)$ in the null space of
$A_0^T$, that is, $A_0^T r = 0$, which is equivalent to the following:
\begin{equation}
    r_1 + \frac{r_2}{1+\alpha_k\tau_e} + \frac{r_3}{1+\alpha_k\tau_i} = 0, \quad k=1,2,3.
\end{equation}
This is a quadratic equation in $\alpha_k$:
\begin{equation}
\begin{aligned}
    f(\alpha_k) \coloneqq{}& r_1 \tau_e \tau_i \alpha_k^2 \\
    &+ \bigl(r_1(\tau_e+\tau_i) + r_2 \tau_i + r_3 \tau_e\bigr) \alpha_k \\
    &+ (r_1 + r_2 + r_3) = 0.
\end{aligned}
\end{equation}
Since a quadratic equation has at most two distinct roots, in order for the equation to hold
for three distinct $\alpha_k$, we must have $f(\alpha_k) \equiv 0$ and $r_1 = r_2 = r_3 = 0$ follows. 
Therefore, the null space of $A_0^T$ is trivial, which implies that $A_0$ has full rank. 
This completes the proof.
\end{proof}

\section{Numerical Results and Discussion}
\label{sec:numerical-results}

We use $C=1\,\mu\mathrm{F}/\mathrm{cm}^2$,
$g_l=0.3\,\mathrm{mS}/\mathrm{cm}^2$, and
$(\bar g_e,\bar g_i)=(8,36)\,\mathrm{mS}/\mathrm{cm}^2$.
The reversal potentials are $(E_e,E_i)=(120,-12)\,\mathrm{mV}$,
with $E_l$ chosen to place rest at the origin.
The activation parameters are $(\theta_e,\theta_i)=(17,10)\,\mathrm{mV}$
and $\kappa=0.5\,\mathrm{mV}^{-1}$; the relaxation times are
$(\tau_e,\tau_i)=(0.4,18)\,\mathrm{ms}$.
As in Section~\ref{sec:structure}, a released state is classified by
whether its zero-input voltage reaches $75\,\mathrm{mV}$ within
$50\,\mathrm{ms}$.

\subsection{Numerical solution}

Numerical algorithms are not the focus of this conference paper;
we outline the scheme. We sample the numerical spiking boundary
more finely near low-cost preparations. At each target $x_f$,
\eqref{eq:three-exponential-terminal-map} determines the coefficients
from the rates, leaving three log-rates to optimize. We recover the
memories and current from the voltage.

We use $0.01\leq\alpha_k\leq10\,\mathrm{ms}^{-1}$,
$\alpha_{k+1}/\alpha_k\geq1.2$, and the preparation box
$\mathcal K=[-12,40]^3\,\mathrm{mV}$. Earlier states are constrained to remain
on the nonspiking side of the numerical release boundary. No additional
current-amplitude or charge-balance constraint is imposed.

The largest cumulative work can occur before release. Writing
\begin{equation}
 W(t)=\int_{-\infty}^{t}i(\tau)v(\tau)\,d\tau,
 \qquad J_{\mathrm{pk}}=\sup_{t\leq0}W(t),
\end{equation}
we compare $W(0)$, the remote-past value zero, and the stationary
values detected through $i(t)v(t)=0$. The capacitor work is integrated
analytically; the ionic work is evaluated by quadrature over the
infinite past.

At a fixed target $x_f$, write $x_\alpha$ and $W_\alpha$ for the
reconstructed state and cumulative work. Introducing a scalar upper
bound $\eta$ gives the epigraph problem
\begin{equation}
\label{eq:numerical-epigraph}
\begin{aligned}
\inf_{\alpha,\eta}\quad &\eta\\
\text{subject to:}\quad
&W_\alpha(t)\leq\eta,\quad t\leq0,\\
&x_\alpha(t)\in\mathcal K,\quad t\leq0,\\
&g_{\mathrm{rel}}(x_\alpha(t))\leq0,\quad t<0,\\
&\alpha\in\mathcal A.
\end{aligned}
\end{equation}
Here $\mathcal A$ contains the rate triples satisfying the stated
bounds and separation, and $g_{\mathrm{rel}}\leq0$ denotes the
nonspiking side of the computed release boundary, including its
numerical threshold bracket. Rest and $x_\alpha(0)=x_f$ are already
enforced by the exponential parameterization. For any admissible $\alpha$,
the smallest feasible $\eta$ is exactly $J_{\mathrm{pk}}$.

Numerically, we approximate the all-time inequalities using detected
work extrema and refined state samples. A quadratic subproblem uses
their linearizations to compute a first-order step.
After restoring feasibility, backtracking checks actual
cost and constraints before accepting a decreasing step. We recompute
work extrema at each trial and add time samples near active
constraints. We compare costs across targets and several rate
initializations per target. Independent work quadrature and zero-input
releases from sampled preparation states check the selected histories.

\subsection{The computed preparation}

\begin{figure*}[t]
    \centering
    \begin{minipage}[t]{0.48\textwidth}
        \centering
        \begin{minipage}[c][0.90\linewidth][c]{\linewidth}
            \centering
            \includegraphics[width=\linewidth]{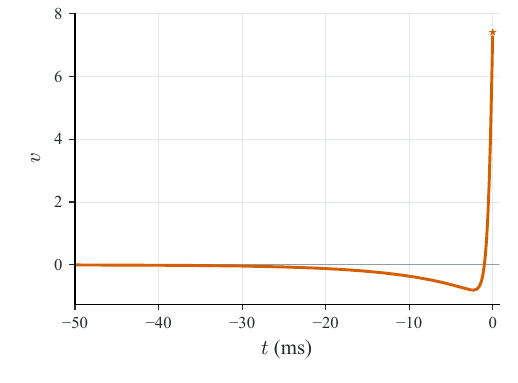}
        \end{minipage}
        \caption{Three-exponential voltage history. The displayed past is
        truncated; work and terminal states include the infinite past.
        Voltage in mV; time in ms.}
        \label{fig:optimal-voltage}
    \end{minipage}\hfill
    \begin{minipage}[t]{0.48\textwidth}
        \centering
        \begin{minipage}[c][0.90\linewidth][c]{\linewidth}
            \centering
            \includegraphics[width=0.90\linewidth]{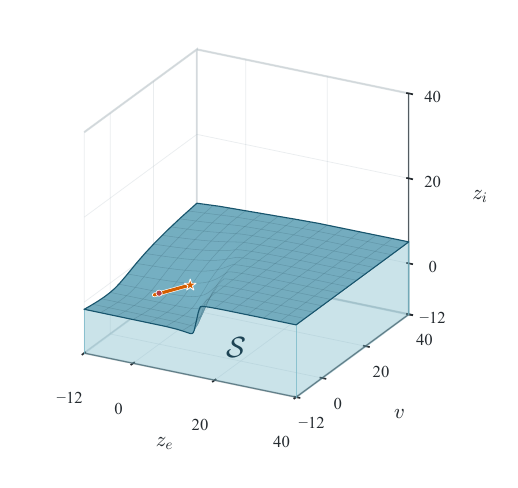}
        \end{minipage}
        \caption{The same preparation (orange) on the spiking region of
        Fig.~\ref{fig:spiking-manifold}. The red dot marks rest; the star
        marks the terminal numerical threshold bracket. Coordinates are in mV.}
        \label{fig:optimal-trajectory}
    \end{minipage}
\end{figure*}

Figure~\ref{fig:optimal-voltage} shows the least costly preparation
found with three exponentials. The voltage first falls to about
$-0.8\,\mathrm{mV}$ and then rises rapidly.
The initial hyperpolarization lowers both memories; during the
final rise, the fast memory recovers while the slow memory remains
negative. Excitatory conductance therefore increases while inhibitory
conductance stays below its resting value.

Figure~\ref{fig:optimal-trajectory} shows the corresponding path to the
spiking boundary, with inhibitory memory still below rest at contact.
We obtain
\begin{equation}
\label{eq:numerical-best-preparation}
\begin{aligned}
    x_f&\simeq(7.4,\,3.7,\,-0.146)^{\top}\,\mathrm{mV},\\
    J_{\mathrm{pk}}&\simeq34.95\,\mathrm{pJ}/\mathrm{cm}^2.
\end{aligned}
\end{equation}
Release depends on all three terminal coordinates, so voltage alone
does not specify the threshold.

\begin{figure*}[t]
    \centering
    \begin{minipage}[c]{0.48\textwidth}
        \centering
        \includegraphics[width=\linewidth,trim=0 20bp 0 0,clip]{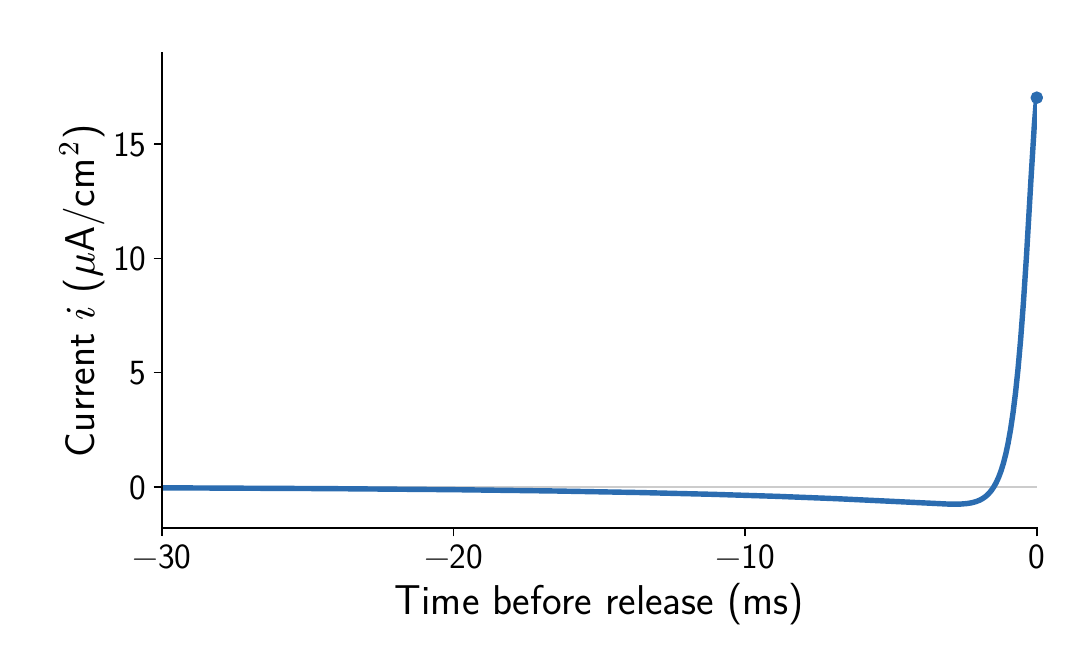}\\[-2pt]
        {\small (a) Reconstructed EI input current}
    \end{minipage}\hfill
    \begin{minipage}[c]{0.48\textwidth}
        \centering
        \includegraphics[width=\linewidth,trim=0 45bp 0 0,clip]{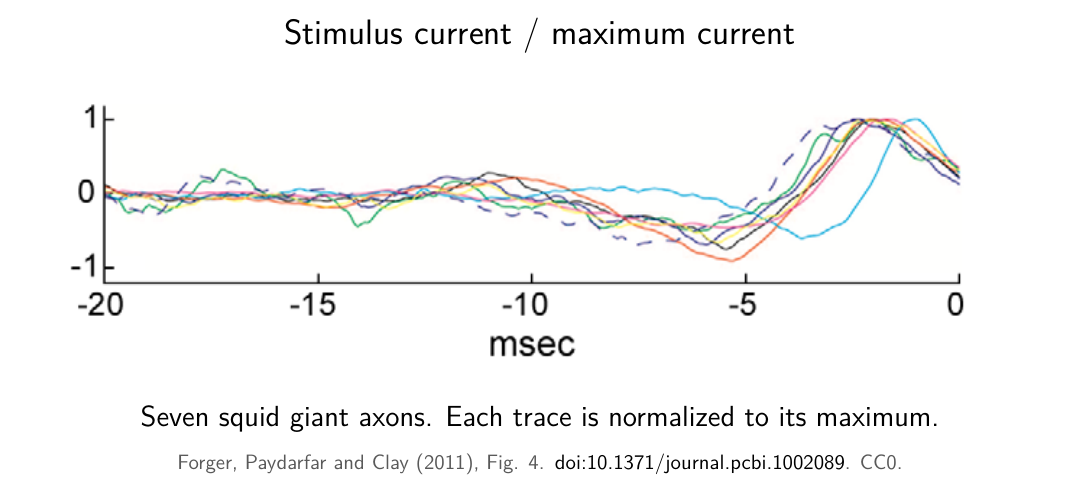}\\[-2pt]
        {\small (b) Experimental stimulus currents}
    \end{minipage}
    \caption{Comparison of applied currents. (a) Current reconstructed
    from the three-exponential preparation, with release at $t=0$.
    (b) Spike-triggered averages
    from seven squid giant axons, reproduced from
    \cite[Fig.~4]{ForgerPaydarfarClay2011}.}
    \label{fig:current-comparison}
\end{figure*}

For this preparation, the largest cumulative work is attained at
release, so $J_{\mathrm{pk}}=W(0)$.
The endpoint represents threshold contact to numerical resolution.
Nearby states with slightly lower and higher inhibitory memory have
spiking and nonspiking zero-input releases, respectively.

For each admissible voltage waveform, our variational formulation concerns 
two outcomes: its associated energy requirement, a continuous quantity, and whether it elicits a 
spike, a binary outcome. In Section~\ref{sec:real-exponential-approximation}, 
we proved that three exponentials suffice to reach any terminal state, 
covering all states in the spiking region. This reachability result does not,
however, establish that three exponentials suffice to attain the optimal energy value.

We therefore increase the number of exponentials to examine whether the 
additional freedom in the voltage waveform can further reduce the energy 
requirement. We repeat the optimization with $N=4,5,6$, under the same 
constraints. We use both independent random initializations and 
warm starts from the solutions obtained with fewer exponentials. 
As shown in Table~\ref{tab:peak-work-vs-N}, increasing $N$ from 
$3$ to $6$ reduces the lowest energy cost by less than $0.001\%$ and the 
waveform retains its hyperpolarizing and depolarizing phases. Although these
results do not establish global optimality, they suggest that three 
exponentials already capture the essential features of low-cost preparations.

\begin{table}[t]
    \centering
    \caption{Smallest peak work found ($\mathrm{pJ}/\mathrm{cm}^2$).}
    \label{tab:peak-work-vs-N}
    \begin{tabular}{crr}
        \hline
        $N$ & Independent starts & Including lower-order starts\\
        \hline
        3 & 34.954505 & 34.954505\\
        4 & 34.954430 & 34.954430\\
        5 & 34.958128 & 34.954421\\
        6 & 34.958052 & 34.954418\\
        \hline
    \end{tabular}
\end{table}

Figure~\ref{fig:current-comparison}(a) shows the input current reconstructed
from the optimized voltage history. It exhibits a negative phase followed by
a positive phase. This is the main significance of our result. First,
it matches the experimentally determined optimal stimulus currents from seven 
squid giant axons in Figure~\ref{fig:current-comparison}(b)~\cite{ForgerPaydarfarClay2011}, where the authors assess efficiency 
through root mean square current. Each experimental trace is normalized to its 
own maximum. In comparison, the computed current has a weaker negative phase 
relative to the positive phase. Second, it aligns with our earlier findings 
in~\cite{SepulchreTong2025} and seems to be a direct consequence of the mixed 
feedback architecture of excitable systems.

\section{Conclusion}
We formulated the problem of finding the minimum energy required to trigger 
a spike, in the spirit of the classical notion of required supply in 
dissipativity theory. We studied this optimal control
problem in a three-dimensional excitable system. The numerical solution yields 
an optimum consistent with experimentally reported results. This is an 
encouraging step toward a variational principle of excitability.
\section*{Acknowledgment}
OpenAI Codex was used for grammar checking and assistance with figure
plotting in Sections~\ref{sec:structure}
and~\ref{sec:numerical-results}.

\bibliographystyle{IEEEtran}
\bibliography{references}

\end{document}